\documentclass[12pt,a4paper,twoside]{article}

\usepackage{amsmath,amsfonts,amsthm,amsopn,color,amssymb}
\usepackage{palatino}
\usepackage{graphicx}

\usepackage[ps2pdf,colorlinks=true,urlcolor=blue,
citecolor=red,linkcolor=blue,linktocpage,pdfpagelabels,bookmarksnumbered,bookmarksopen]{hyperref}

\newcommand{\C}{\mathbb{C}}
\newcommand{\R}{\mathbb{R}}

\newcommand{\RN}{{\mathbb{R}^N}}

\newcommand{\de}{\partial}

\renewcommand{\le}{\leqslant}
\renewcommand{\ge}{\geqslant}

\renewcommand{\d }{\delta }
\newcommand{\e}{{\rm e}}
\newcommand{\vfi}{\varphi}
\newcommand{\g }{\gamma }

\newcommand{\n }{\nabla }

\renewcommand{\t}{\theta}
\renewcommand{\O}{\Omega}

\newcommand{\G}{\Gamma}

\newcommand{\A}{{\cal A}}
\newcommand{\B}{{\cal B}}
\renewcommand{\C}{{\cal C}}

\renewcommand{\H}{H^1(\RN)}
\newcommand{\HH}{\mathbb{H}}

\newcommand{\Ne}{\mathcal{N}}

\newcommand{\SU}{\sum_{i=1}^3}
\newcommand{\N}{\mathbb{N}}

\renewcommand{\o}{\omega}

\newcommand{\bu}{{\bf u}}
\newcommand{\bun}{{\bf u}_n}

\newcommand{\bvn}{{\bf v}_n}
\newcommand{\bwn}{{\bf w}_n}
\newcommand{\bzn}{{\bf z}_n}
\newcommand{\uuu}{(u_1,u_2,u_3)}

\newcommand{\uuun}{(u_{1,n},u_{2,n},u_{3,n})}
\newcommand{\uun}{u_{1,n} u_{2,n} u_{3,n}}

\newcommand{\irn }{\int_{\RN}}

\def\bbm[#1]{\mbox{\boldmath $#1$}}

\newtheorem{theorem}{Theorem}[section]
\newtheorem{lemma}[theorem]{Lemma}

\newtheorem{definition}[theorem]{Definition}

\renewenvironment{proof}{\noindent{\textbf{Proof\quad}}}{$\hfill\square$\vspace{0.2 cm}\\}
\newenvironment{proofmain}{\noindent{\textbf{Proof of Theorem  \ref{main}\quad}}}{$\hfill\square$\vspace{0.2 cm}\\}
\newenvironment{proofmain2}{\noindent{\textbf{Proof of Theorem \ref{main2}\quad}}}{$\hfill\square$\vspace{0.2 cm}\\}
\newenvironment{proofmain3}{\noindent{\textbf{Proof of Theorem  \ref{main3}\quad}}}{$\hfill\square$\vspace{0.2 cm}\\}
\newenvironment{proofmain4}{\noindent{\textbf{Proof of Theorem \ref{main4}\quad}}}{$\hfill\square$\vspace{0.2 cm}\\}

\title{{\sc Ground states for a system of nonlinear Schr\"odinger equations with three waves interaction\footnote{The authors are supported by M.I.U.R. -
P.R.I.N. ``Metodi variazionali e topologici nello studio di
fenomeni non lineari''}}}

\author{A. Pomponio\thanks{Dipartimento di Matematica, Politecnico di
Bari, Via E. Orabona 4, I-70125 Bari, Italy, e-mail: {\tt
a.pomponio@poliba.it}}}

\date{}

\begin{document}

\maketitle

\begin{abstract}
We consider a system of nonlinear Schr\"odinger equations with three waves interaction studying the existence of ground state solutions. In particular, we find a {\it vector} ground state, namely a ground state $\uuu$ such that $u_i\neq 0$, for all $i=1,2,3$.
\end{abstract}

\section{Introduction}

When light is scattered from an atom or molecule, most photons are elastically scattered (Rayleigh scattering), such that the scattered photons have the same energy (frequency) and wavelength as the incident photons. A small fraction of the scattered light (approximately 1 in 10 million photons), however,  is scattered by an excitation, with the scattered photons having a frequency different from, and usually lower than, the frequency of the incident photons. This process is know as Raman effect or Stimulated Raman Scattering. The Raman amplification is based on this phenomenon: when a lower frequency signal photon induces the inelastic scattering of a higher-frequency pump photon in an optical medium in the nonlinear regime, as a result of this, another signal photon is produced, with the surplus energy resonantly passed to the vibrational states of the medium. This process, as with other stimulated emission processes, allows all-optical amplification. In addition to applications in nonlinear and ultrafast optics, Raman amplification is used in optical telecommunications, allowing all-band wavelength coverage and in-line distributed signal amplification. For physics and engineering	aspects about this topic see, for example, \cite{HA}. 

This phenomenon, in a suitable mathematical contest, has bee treated in \cite{CC,CC2}, where a new set of equations describing nonlinear interaction between a laser beam and a plasma has been derived. From a physics point of view, when an incident laser field enters a plasma, it is backscattered
by a Raman type process. Then these two waves interact to create an electronic plasma wave. The three waves combine to
create a variation of the density of the ions which has itself an influence on the three proceedings waves. The system
describing this phenomenon is composed by three Schr\"odinger equations coupled to a wave equation.

Later on, this model, slightly modified, has been studied in \cite{CCO1,CCO2}, where in particular the orbital stability of solitary waves is studied for the following system of nonlinear Schr\"odinger equations:
\begin{equation}\label{eq:i}
\left\{
\begin{array}{l}
i\de_t v_1=-\Delta v_1 - |v_1|^{p-1}v_1-\g v_3 \overline{v_2},
\\
i\de_t v_2=-\Delta v_2 - |v_2|^{p-1}v_2-\g v_3 \overline{v_1},
\\
i\de_t v_3=-\Delta v_3 - |v_3|^{p-1}v_3-\g v_1  v_2,
\end{array}
\right.
\end{equation}
where $v_1$, $v_2$ and $v_3$ are complex valued functions of $(t,x)\in \R\times \RN$, $N=1,2,3$, $1<p<1+4/N$, $\g>0$. 
More precisely, if $\vfi\in H^1(\RN)$ is the unique positive radial (least-energy) solution of 
\[
-\Delta v +\o v - |v|^{p-1}v=0 \qquad \hbox{in }\RN,
\]
with $\o>0$ (see \cite{BL1,K}), the authors study the orbital stability of the following standing wave solutions: $(\e^{i\o t}\vfi,0,0)$, $(0,\e^{i\o t}\vfi,0)$ and $(0,0,\e^{i\o t}\vfi)$. They 
show that  $(\e^{i\o t}\vfi,0,0)$ and $(0,\e^{i\o t}\vfi,0)$ are orbitally stable, for any $\g>0$, while 
$(0,0,\e^{i\o t}\vfi)$ is orbitally stable, if $0<\g<\g^*$, and it is orbitally unstable, if $\g>\g^*$, for a suitable positive constant $\g^*=\g^*(N,p,\o)$.

If, moreover, among the solutions of \eqref{eq:i}, we look for standing waves, namely solutions of the type
\[
(v_1(t,x),v_2(t,x),v_3(t,x))=(\e^{i \o_1 t}u_1(x),\e^{i \o_2 t}u_2(x),\e^{i \o_3 t}u_3(x)),
\]
where $u_1$, $u_2$ and $u_3$ are real valued functions of $x\in \RN$, it is easy to see that, if $\o_3=\o_1+\o_2$, $(\e^{i \o_1 t}u_1(x),\e^{i \o_2 t}u_2(x),\e^{i \o_3 t}u_3(x))$ is a solution of \eqref{eq:i} if $\uuu$ is a solution of the following system:
\begin{equation}\label{eq0}
\left\{
\begin{array}{l}
-\Delta u_1 +\o_1 u_1 - |u_1|^{p-1}u_1=\g u_2 u_3,
\\
-\Delta u_2 +\o_2 u_2 - |u_2|^{p-1}u_2=\g u_1 u_3,
\\
-\Delta u_3 +\o_3 u_3 - |u_3|^{p-1}u_3=\g u_1 u_2.
\end{array}
\right.
\end{equation}
Therefore, the orbital stability of such type of solutions is strictly linked with the existence of a least energy solution of \eqref{eq0}.

Motivated by these previous results and by the above considerations, in this paper we are interested in a generalization of system \eqref{eq0}. More precisely, we consider the following system:
\begin{equation}\label{eq}\tag{${\cal P}_\o$}
\left\{
\begin{array}{ll}
-\Delta u_1 +\o_1 u_1 - |u_1|^{p-1}u_1=\g u_2 u_3	& \hbox{in }\RN,
\\
-\Delta u_2 +\o_2 u_2 - |u_2|^{p-1}u_2=\g u_1 u_3	& \hbox{in }\RN,
\\
-\Delta u_3 +\o_3 u_3 - |u_3|^{p-1}u_3=\g u_1 u_2	& \hbox{in }\RN,
\\
u_1,u_2,u_3\in \H,
\end{array}
\right.
\end{equation}
where $\o_i>0$, $i=1,2,3$, $\g\in \R$ and $2<p<(N+2)/(N-2)$. It is easy to observe that \eqref{eq} possesses {\it trivial} solutions different from $(0,0,0)$. Indeed, by \cite{BL1} there exists a (least--energy) solution $u_i \in \H$ for the single Schr\"odinger equation
\begin{equation*}	
-\Delta u+ \o_i u=|u|^{p-1} u \qquad \hbox{in }\RN.
\end{equation*}
It can be checked immediately that $(u_1,0,0)$, $(0,u_2,0)$ and $(0,0,u_3)$ are non-trivial solutions of \eqref{eq}. We will refer to these solutions as {\it scalar solution} and, of course, they are meaningless. We are looking for other type of solutions that we will call {\it vector solution}. 

We set $\bu=\uuu$ and $\HH=\H\times\H\times \H$.

To fix terminology, we introduce the following definition.
\begin{definition}
A solution of \eqref{eq}, $\bu\in \HH$, $\bu\neq(0,0,0)$ will be called {\it scalar solution} if there exist $i,j=1,2,3$, $i\neq j$, such that $u_i\equiv u_j \equiv 0$; while a solution $\bu\in \HH$ of \eqref{eq} will be called {\it vector solution} if $u_i\neq0$, for $i=1,2,3$.
\end{definition}

System \eqref{eq} has a variational structure and so its solutions can be found as critical points of the functional $I\colon \HH \to \R$ defined as follows:
\begin{equation}\label{eq:I}
I(\bu)=\sum_{i=1}^3 I_i(u_i)-\g \irn u_1 u_2 u_3,
\end{equation}
where, for $i=1,2,3$,
\[
I_i(u_i)=\frac 12 \irn |\n u_i|^2+\o_i u_i^2 
-\frac{1}{p+1}\irn u_i^{p+1}.
\]

We are interested in the existence of a \emph{ground state solution}, namely $\bu\in \HH$,  $\bu\neq (0,0,0)$ which solves
(\ref{eq}) and minimizes the functional $I$ among all possible nontrivial solutions.

As a first step, we prove that, for any $\g\in \R$, the problem \eqref{eq} admits a ground state.

\begin{theorem}\label{main}
For any $\o_i>0$, $i=1,2,3$ and for any $\g\in \R$, there exists a ground state solution $\bu\in \HH$ of \eqref{eq}. Moreover $u_i$ is a radially symmetric function (up to translation), for $i=1,2,3$. 
\end{theorem}

Then we prove that vector solutions exist whenever the absolute value of the coupling parameter $\g$ is sufficiently large.
\begin{theorem}\label{main2}
For any $\o_i>0$, $i=1,2,3$, there exists $\g_0>0$ such that, for any $\g \in \R$ with $|\g|>\g_0$, \eqref{eq} possesses a vector solution $\bu\in \HH$, which is a ground state solution. 
\end{theorem}

The main result of the first part of this paper is Theorem \ref{main2}: up to our knowledge, indeed, this is the first vector solution existence result for problem \eqref{eq}. Theorem \ref{main}, instead, is essentially already known. The symmetry result is due to \cite{BJM} while the existence result is proved in \cite{BL}. For the reader's sake, here we give a different proof of the existence of a ground state and our arguments are based on the constrained minimization over the Nehari manifold. Theorems \ref{main} and \ref{main2} will be proved in Section \ref{se:w}.

By a mathematical point of view, system \eqref{eq} is related with weakly coupled nonlinear Schr\"odinger systems (see for example \cite{PS} and the bibliography therein). However, the peculiarity of \eqref{eq} is due to the presence, in the functional $I$ defined in \eqref{eq:I}, of an integral term without a pre-assigned sign and this, of course, generates several difficulties.

In the second part of the paper, we try to generalize further on system \eqref{eq}. We consider, indeed, the case when we substitute in \eqref{eq} the positive constants $\o_i$ with non-constant positive potentials $V_i(x)$. More precisely, we consider
\begin{equation}\label{eqv}\tag{${\cal P}_{ V}$}
\left\{
\begin{array}{ll}
-\Delta u_1 +V_1(x) u_1 - |u_1|^{p-1}u_1=\g u_2 u_3	& \hbox{in }\RN,
\\
-\Delta u_2 +V_2(x) u_2 - |u_2|^{p-1}u_2=\g u_1 u_3	& \hbox{in }\RN,
\\
-\Delta u_3 +V_3(x) u_3 - |u_3|^{p-1}u_3=\g u_1 u_2	& \hbox{in }\RN,
\\
u_1,u_2,u_3\in \H,
\end{array}
\right.
\end{equation}
where $\g\in \R$, $2<p<(N+2)/(N-2)$ and we assume the following hypotheses on ${\bf V}=(V_1,V_2,V_3)$:
\begin{itemize}
\item[({\bf V1})] for all $i=1,2,3$, $V_i\colon\RN \to\R$ is a measurable function;
\item[({\bf V2})] for all $i=1,2,3$, $V_{i,\infty}:=\lim_{|y|\to\infty}V_i(y) \ge V_i(x)$, for almost every $x\in \RN$, and the inequality is strict in a non-zero measure domain;
\item[({\bf V3})] for all $i=1,2,3$, $0<C_i\le V_i(x)$, for all $x\in \RN$.
\end{itemize}
Such type of assumptions has been introduced in \cite{R} for the study of a single nonlinear Schr\"odinger equation.

With concentration-compactness arguments, we prove the following:
\begin{theorem}\label{main3}
Suppose that ${\bf V}=(V_1,V_2,V_3)$ satisfies ({\bf V1-3}), then for any $\g\in \R$, there exists a ground state solution $\bu\in \HH$ of \eqref{eqv}.
\end{theorem}

As in the first part of the paper, we prove that vector solutions exist whenever the absolute value of the coupling parameter $\g$ is sufficiently large. In this case, moreover, we can also slightly weaken the assumptions on the potential ${\bf V}$, in particular, instead of ({\bf V2}), we make the following assumption:
\begin{itemize}
\item[({\bf V2'})] for all $i=1,2,3$, $V_{i,\infty}:=\lim_{|y|\to\infty}V_i(y) \ge V_i(x)$, for almost every $x\in \RN$, and there exists at least one $i=1,2,3$ such that the inequality is strict in a non-zero measure domain.
\end{itemize}
We remark that ({\bf V2}) requires that all the potentials $V_i$ have the same geometrical behavior and in particular all the three potentials have to be non-constant; at contrary, the assumption ({\bf V2'}) is satisfied if, for example, only one of the three potentials is non-constant and with the right geometry, while the other two could be positive constants.

The following theorem holds:
\begin{theorem}\label{main4}
Suppose that ${\bf V}=(V_1,V_2,V_3)$ satisfies ({\bf V1}),({\bf V2'}) and ({\bf V3}), then there exists $\g_0>0$ such that, for any $\g\in \R$ with $|\g|>\g_0$, \eqref{eqv} possesses a vector solution $\bu\in \HH$, which is a ground state solution. 
\end{theorem}

Up to our knowledge, these last two theorems are the first results for a system of nonlinear Schr\"odinger equations with three waves interaction in presence of non-constant potentials. Theorems \ref{main3} and \ref{main4} will be proved in Section \ref{se:V}.

\medskip

\begin{center}{\bf Notation}\end{center}
\begin{itemize} 
\item If $r>0$ and $x_0 \in \RN$, $B_r (x_0):= \left\{ x\in\RN : |x- x_0| <r \right\}$.
We denote with $B_r$ the ball of radius $r$ centered in the origin.
\item We set $\bu=\uuu$, $\bun =\uuun$ and ${\bf V}=(V_1,V_2,V_3)$.
\item We denote by $\|\cdot \|$ the standard norm of $\H$.
\item We set $\HH=\H \times \H \times \H$ and, for any $\bu\in \HH$, we set $\|\bu\|^2=\SU \|u_i\|^2$.
\item For any $1\le s\le +\infty$, we denote by $\|\cdot\|_s$ the usual norm of the Lebesgue space $L^s(\RN)$.
\item By $C_i$ and $c_i$, we denote generic positive constants, which may also vary from line to line.
\item By $o_n(1)$ we denote a quantity which vanishes to zero as $n\to +\infty$. 
\end{itemize}

\section{The constant potential case}\label{se:w}

We prove Theorem \ref{main} by a constrained minimization over the Nehari manifold.
Let us define $G:\HH\to \R$ such that
\begin{align*}
G (\bu) &= I'(\bu)[\bu]
= \sum_{i=1}^3 I'_i(u_i)[u_i]-3\g \irn u_1 u_2 u_3
\\
&=\sum_{i=1}^3 \irn |\n u_i|^2+\o_i u_i^2 - |u_i|^{p+1}-3\g \irn u_1 u_2 u_3,
\end{align*}
then, any critical point, $\bu$, of $I$ satisfies the following equality:
\[
G (\bu)=0.
\]
We denote by $\Ne$ the so called Nehari manifold of $I$, namely
\[
\Ne :=\left\{ \bu\in \HH \setminus \{(0,0,0)\} \;\Big{|}\; G (\bu)=0 \right\}.
\]

The next lemmas give useful informations on the Nehari manifold.
\begin{lemma}\label{le:N1}
For any $\g\in \R$, there exists a positive constant $C_\g$, such that for all $\bu \in \Ne$, $\|\bu\|\ge C_\g$.
\end{lemma}

\begin{proof}
The conclusion follows immediately, observing that, for any $\bu\in \Ne$, we have
\begin{align*}
\|\bu\|^2&=\SU \|u_i\|^2 \le C\SU \irn |\n u_i|^2 +\o_i u_i^2
\\
&= C \left(\SU \irn |u_i|^{p+1} - 3\g \irn u_1 u_2 u_3\right)
\\
&\le C \left(\SU \|u_i\|_{p+1}^{p+1} +3 |\g|\cdot\|u_1\|_3 \|u_2\|_3 \|u_3\|_3\right)
\\
&\le C \left(\SU \|u_i\|^{p+1} +3 |\g|\cdot\|u_1\| \|u_2\| \|u_3\|\right).
\end{align*}
\end{proof}

\begin{lemma}\label{le:N2}
$\Ne$ is a $C^1$ manifold and it is a natural constraint for the functional $I$, namely each critical point of $I_{|\Ne}$ is a critical point of the unconstrained functional~$I$.
\end{lemma}

\begin{proof}
Let $\bu\in \Ne$, we have:
\begin{align*}
G'(\bu)[\bu]&=
\SU \irn 2|\n u_i|^2+2\o_i u_i^2 - (p+1)|u_i|^{p+1}
-9\g \irn u_1 u_2 u_3
\\
&=\SU \irn -|\n u_i|^2-\o_i u_i^2 - (p-2)|u_i|^{p+1}
<-C\|\bu\|^2<-C_\g<0.
\end{align*}
Hence we get the conclusion.
\end{proof}

\begin{lemma}\label{le:N3}
For any $\bu\in\HH$, $\bu\neq (0,0,0)$ there exists a unique number $\bar t>0$ such that $\bar t \bu\in \Ne$ and
\[
I(\bar t \bu)=\max_{t \ge 0}I(t \bu).
\]
\end{lemma}

\begin{proof}
Let $\bu\in\HH$, $\bu\neq (0,0,0)$. For any $t>0$, we set
\[
f(t)=I(t\bu)
=\frac{\A}{2}t^2-\frac{\B}{p+1}t^{p+1}-\C t^3,
\]
where
\begin{align*}
\A&=\SU \irn |\n u_i|^2+\o_i u_i^2>0,
\\
\B&=\SU \irn |u_i|^{p+1}>0,
\\
\C&=\g \irn u_1 u_2 u_3 \in \R.
\end{align*}
We have to show that $f:\R_+ \to \R$ admits a unique maximum. 
Let us observe that
\[
\lim_{t\to 0^+}\frac{f(t)}{t^2}= \frac{\A}{2}>0,\qquad \lim_{t\to +\infty} f(t)=-\infty,
\]
so there exists at least a $\bar t>0$ such that 
\[
f(\bar t)=\max_{t\ge 0} f(t).
\]
Since any maximum point of $f$ satisfies
\begin{equation}\label{eq:f}
f'(t)=\A t-\B t^{p}-3 \C t^2=0 \Longleftrightarrow
\A =\B t^{p-1}+3 \C t,
\end{equation}
we conclude if we prove that $\bar t$ is the unique solution of \eqref{eq:f}. 
\\
We have to distinguish two cases. If $\C\ge 0$, then it is easy to see that \eqref{eq:f} admits a unique solution. Suppose, instead, that $\C<0$ and set ${\cal D}=\!-\C>0$. 
\\
For $t>0$, we set 
\[
g(t)=\frac{f'(t)}{t}= \A -\B t^{p-1}+3 {\cal D} t.
\]
Let $t_0=\min\{t> 0\mid g(t)=0\}$ and suppose, by contradiction, that there exists $t_1>t_0$ such that
\begin{equation}\label{eq:t01}
g(t_0)=g(t_1)=0.
\end{equation}
Since $\lim_{t\to 0} g(t)=\A>0$, then $g'(t_0)\le 0$ and so
\begin{equation}\label{eq:t0}
3 {\cal D}\le (p-1)\B t_0^{p-2}.
\end{equation}
Moreover by \eqref{eq:t01}, there exists $t_2\in (t_0,t_1)$, such that $g'(t_2)=0$. Hence
\[
3 {\cal D}=(p-1)\B t_2^{p-2}> (p-1)\B t_0^{p-2},
\]
and we get a contradiction with \eqref{eq:t0}.
\end{proof}

According to the definition of \cite{L1}, we say that a sequence $\{\bun\}_n$ vanishes if, for all $r>0$
\[
\lim_{n \to +\infty} \sup_{\xi \in \RN} \int_{B_r(\xi)}u_{1,n}^2+u_{2,n}^2+u_{3,n}^2=0.
\]
\begin{lemma}\label{le:nonvan}
Any bounded sequence $\{\bun\}_n\subset \Ne$ does not vanish.
\end{lemma}

\begin{proof}
Suppose by contradiction that $\{\bun\}_n$ vanishes, then, in particular there exists $\bar r>0$ such that
\[
\lim_{n \to +\infty} \sup_{\xi \in \RN} \int_{B_{\bar r}(\xi)}u_{i,n}^2=0,
\quad \hbox{for all }i=1,2,3.
\]
Then, by \cite[Lemma 1.1]{L2}, we infer that $u_{i,n}\to 0$ in $L^s(\RN)$, $i=1,2,3$, for any $2<s<2^*$. Since $\{\bun\}_n\subset \Ne$, we have that $\bun \to 0$ in $\HH$, contradicting Lemma \ref{le:N1}.
\end{proof}

Now we are ready to prove Theorem \ref{main}.

\begin{proofmain}
It is easy to see that for any $\bu\in \Ne$ we get
\begin{equation}\label{eq:IN}
I (\bu)= \SU  \irn \frac 16\left(|\n u_i|^2+\o_i u_i^2\right) +\frac{p-2}{3(p+1)}|u_i|^{p+1},
\end{equation}
hence, by Lemma \ref{le:N1} and since $p>2$, we infer that
\[
m=\inf_{\bu \in \Ne} I(\bu)>0. 
\]
We have to show that this infimum is achieved as a minimum. 
\\
Let $\{\bun\}_n\subset \Ne$ be a minimizing sequence. By standard arguments (see \cite{Wi}), we can suppose that $\{\bun\}$ is a Palais-Smale sequences for $I$ at level $m$, namely
\begin{align}
I(\bun) \to m, &\qquad\hbox{ as }n \to \infty,\label{eq:PS1}
\\ 
I'(\bun) \to 0, &\qquad\hbox{ as }n \to \infty.\label{eq:PS2}
\end{align}
By \eqref{eq:IN} and \eqref{eq:PS1}, we know that $\{\bun\}_n$ is a bounded sequence in $\HH$. Therefore, by Lemma \ref{le:nonvan}, $\{\bun\}_n$ does not vanish, namely there exist $C,r>0$, $\{\xi_n\}_n\subset \RN$ such that
\begin{equation}\label{eq:nonvan}
\int_{B_r(\xi_n)}u_{1,n}^2+u_{2,n}^2+u_{3,n}^2\ge C, \hbox{ for all }n\ge 1.
\end{equation}
Due to the invariance by translations, without loss of generality, we can assume that $\xi_n=0\in \RN$, for every $n$.
\\
Since $\{\bun\}_n$ is bounded in $\HH$, there exist $u_1,u_2,u_3\in \H$ such that, up to a subsequence, for $i=1,2,3$,
\begin{align*}
&u_{i,n} \rightharpoonup u_i, \quad \hbox{in }\H; 
\\
&u_{i,n} \to u_i, \quad \hbox{a.e. in }\RN; 
\\
&u_{i,n} \to u_i, \quad \hbox{in } L^s_{\mathrm{loc}}(\RN),\; \;1\le s<2^*.
\end{align*}
Hence, by \eqref{eq:PS2}, we infer that $I'(\bu)=0$.
\\
By \eqref{eq:nonvan}, moreover, we can argue that there exists $i=1,2,3$ such that $u_i\neq 0$, namely $\bu\neq (0,0,0)$, and so we can conclude that $\bu\in\Ne$. By the weak lower semicontinuity, we get
\begin{align*}
m &\le I(\bu)
=\SU  \irn \frac 16\left(|\n u_i|^2+\o_i u_i^2\right) +\frac{p-2}{3(p+1)}|u_i|^{p+1}
\\
& \le \liminf_{n \to +\infty} \SU  \irn \frac 16\left(|\n u_{i,n}|^2+\o_i u_{i,n}^2\right) 
+\frac{p-2}{3(p+1)}|u_{i,n}|^{p+1}
\\
&= \liminf_{n \to +\infty} I(\bun)=m,
\end{align*}
hence $\bu$ is a ground state for the problem \eqref{eq}.
\\
Finally, by \cite{BJM}, we infer that $u_i$ is a radially symmetric function (up to translation), for $i=1,2,3$.
\end{proofmain}

Let us now prove the existence of vector ground state. Here we use some ideas of \cite{MMP}.

\begin{proofmain2}
We start proving that there exists $\g_0>0$ such that \eqref{eq} admits a vectorial ground state, for any $\g>\g_0$. 
\\
Let $\bar u_i\in \H$ be the positive radial ground state of 
\[
-\Delta u+\o_i u= |u|^{p-1}u, \quad \hbox{in }\RN, \qquad u \in \H.
\]
We conclude if we show that there exists $\bu\in \Ne$ such that
\begin{equation}\label{eq:vec}
I(\bu)<\min\big\{I(\bar u_1,0,0),I(0,\bar u_2,0),I(0,0,\bar u_3)\big\}
=\min_{i=1,2,3}I_i(\bar u_i).
\end{equation}
By Lemma \ref{le:N3}, we know that there exists $t_\g>0$ such that $t_\g (\bar u_1,\bar u_2,\bar u_3)\in \Ne$, namely
\[
\SU \irn |\n \bar u_i|^2+\o_i \bar u_i^2
=t_\g^{p-1}\SU \irn |\bar u_i|^{p+1}
+3 t_\g \g \irn \bar u_1 \bar u_2 \bar u_3.
\]
Since the last integral is strictly positive, it is easy to see that
\begin{equation}\label{eq:tg+}
\lim_{\g\to +\infty} t_\g=0.
\end{equation}
Moreover since
\[
I (t_\g (\bar u_1,\bar u_2,\bar u_3))
= \SU   \irn \frac{t_\g^2}6\left(|\n \bar u_i|^2+\o_i \bar u_i^2\right) +\frac{(p-2)t_\g^{p+1}}{3(p+1)}|\bar u_i|^{p+1},
\]
by \eqref{eq:tg+}, for $\g$ positive and sufficiently large, \eqref{eq:vec} is satisfied.
\\
Let us now prove the same conclusion also for $\g$ negative and sufficiently large in modulus. 
\\
By Lemma \ref{le:N3}, we know that there exists $t_\g>0$ such that $t_\g (-\bar u_1,\bar u_2,\bar u_3)\in \Ne$, namely
\[
\SU \irn |\n \bar u_i|^2+\o_i \bar u_i^2
=t_\g^{p-1}\SU \irn |\bar u_i|^{p+1}
-3 t_\g \g \irn \bar u_1 \bar u_2 \bar u_3.
\]
Observing that
\begin{equation*}
\lim_{\g\to -\infty} t_\g=0,
\end{equation*}
we conclude arguing as in the previous case.
\end{proofmain2}

\section{The non-constant potential case}\label{se:V}

In this section, in not stated differently, we will always assume that ${\bf V}$ assume ({\bf V1-3}).

Solutions of \eqref{eqv} are critical points of the functional $I_{\bf V}\colon \HH \to \R$ so defined:
\[
I_{\bf V}(\bu)=\sum_{i=1}^3 I_{V_i}(u_i)-\g \irn u_1 u_2 u_3,
\]
where, for $i=1,2,3$,
\[
I_{V_i}(u_i)=\frac 12 \irn |\n u_i|^2+V_i(x) u_i^2 
-\frac{1}{p+1}\irn u_i^{p+1}.
\]
As in the first part of the paper, let us define 
\begin{align*}
G_{\bf V} (\bu) &= I'_{\bf V}(\bu)[\bu]
= \sum_{i=1}^3 I'_i(u_i)[u_i]-3\g \irn u_1 u_2 u_3
\\
&=\sum_{i=1}^3 \irn |\n u_i|^2+V_i(x) u_i^2 - |u_i|^{p+1}-3\g \irn u_1 u_2 u_3,
\end{align*}
then, any critical point, $\bu$, of $I_{\bf V}$ satisfies the following equality:
\[
G_{\bf V} (\bu)=0.
\]
We denote by $\Ne_{\bf V}$ the Nehari manifold of $I_{\bf V}$, namely
\[
\Ne_{\bf V} :=\left\{ \bu\in \HH \setminus \{(0,0,0)\} \;\Big{|}\;  G_{\bf V} (\bu)=0 \right\}.
\]

The following lemmas describe some properties of the Nehari
manifold $\Ne_{\bf V}$. The proofs are very similar to those of the constant potential case and we omit them.
\begin{lemma}\label{le:NV1}
For any $\g\in \R$, there exists a positive constant $C_\g$, such that for all $\bu \in \Ne_{\bf V}$, $\|\bu\|\ge C_\g$.
\end{lemma}

\begin{lemma}\label{le:NV2}
$\Ne_{\bf V}$ is a $C^1$ manifold and it is a natural constraint for the functional $I_{\bf V}$.
\end{lemma}

\begin{lemma}\label{le:NV3}
For any $\bu\in\HH$, $\bu\neq (0,0,0)$ there exists a unique number $\bar t>0$ such that $\bar t \bu\in \Ne_{\bf V}$ and
\[
I_{\bf V}(\bar t \bu)=\max_{t \ge 0} I_{\bf V}(t \bu).
\]
\end{lemma}

For all ${\bf V}$, we assume the following definition
\begin{equation*}
c_{\bf V}:=\inf_{\bu\in\Ne_{\bf V}} I_{\bf V}(\bu),
\end{equation*}
so that our goal is to find $\bar \bu\in\Ne_{\bf V}$ such that $ I_{\bf V}(\bar \bu)=c_{\bf V}$, from which we would deduce that $\bar \bu$ is a ground state solution of \eqref{eqv}.

Let us recall some preliminary lemmas which can be obtained by
using the same arguments as in \cite{R} (see also \cite{AP}).

As a consequence of the Lemma \ref{le:NV3}, we are allowed to
define the map $t:\HH\setminus\{(0,0,0)\}\to\R_+$ such that for any
$\bu\in\HH,$ $\bu\neq (0,0,0):$
\begin{equation*}
I_{\bf V}\big( t_{\bu} \bu\big)=\max_{t\ge 0}  I_{\bf V}(t \bu).
\end{equation*}
\begin{lemma}\label{le:ccc2}
The following equalities hold
\begin{equation*}
c_{\bf V}=\inf_{g\in \G} \max_{t\in [0,1]}
I_{\bf V}(g(t))=\inf_{\bu\neq (0,0,0)} \max_{t \ge 0} I_{\bf V}( t\bu),
\end{equation*}
where 
\begin{equation*}
\G=\left\{g\in C\big([0,1],\HH\big) \mid g(0)=(0,0,0),\;I_{\bf V}(g(1))\le 0,
\;g(1)\neq (0,0,0)\right\}.
\end{equation*}
\end{lemma}

\begin{lemma}\label{le:>d}
Let $\{\bun\}_n \subset \HH$,  such that $\|\bun\|\ge C>0$, for all $n\ge 1,$ and
\begin{equation*}
\max_{t\ge 0} I_{\bf V}(t \bun) \le c_{\bf V}+\d_n,
\end{equation*}
with $\d_n\to0^+.$ Then, there exist a sequence $\{y_n\}_n \subset \RN$ and two positive numbers $R,\;\mu>0$ such that
\begin{equation*}
\liminf_n \int_{B_R(y_n)}u_{1,n}^2+u_{2,n}^2+u_{3,n}^2 \, d x >\mu.
\end{equation*}
\end{lemma}

\begin{lemma}\label{le:t-bdd}
Let $\{\bun\}_n  \subset \HH$ and  $\{t_n\}_n\subset \R_+$ such that $0<C_1\le \|\bun\|\le C_2$ and
\begin{equation*}
I_{\bf V}(t_n\bun)= \max_{t\ge 0}I_{\bf V}( t \bun) \to c_{\bf V},\; \hbox{ as }n
\to \infty.
\end{equation*}
Then the sequence $\{t_n\}_n$ possesses a bounded
subsequence in $\R$.
\end{lemma}

\begin{proof}
We have
\begin{align*}
C&\ge \SU \irn |\n u_{i,n}|^2 + V_i(x)u_{i,n}^2 
\\
&=t_n \left(t_n^{p-2} \SU \irn |u_{i,n}|^{p+1}+3\g \irn \uun \right).
\end{align*}
The conclusion follows from Lemma \ref{le:>d}.
\end{proof}

\begin{lemma}\label{le:vn}
Suppose that $V_i,$ $V_{i,n} \in L^\infty(\RN)$, for $i=1,2,3$ and for all $n\ge 1$.
\\
If $V_{i,n}\to V_i$ in $L^\infty(\RN)$ then $c_{{\bf V}_n}\to c_{\bf V}$, where ${\bf V}=(V_1,V_2,V_3)$ and ${\bf V}_n=(V_{1,n},V_{2,n},V_{3,n})$.
\end{lemma}

Let $I_\infty\colon \HH \to \R$ be the functional so defined:
\[
I_\infty (\bu)=\sum_{i=1}^3 I_{i,\infty}(u_i)-\g \irn u_1 u_2 u_3,
\]
where, for $i=1,2,3$,
\[
I_{i,\infty}(u_i)=\frac 12 \irn |\n u_i|^2+V_{i,\infty}(x) u_i^2 
-\frac{1}{p+1}\irn u_i^{p+1}.
\]
Moreover, if we set ${\bf V}_\infty=(V_{1,\infty},V_{2,\infty},V_{3,\infty})$, we denote
\[
c_\infty :=c_{{\bf V}_\infty}.
\]

As in \cite{R}, we have
\begin{lemma}\label{le:cinfty}
If ${\bf V}$ satisfies ({\bf V1-3}), then, for any $\g\in \R$, we get $c_{\bf V}<c_\infty$.
\end{lemma}

\begin{proof}
By Theorem \ref{main}, there exists $\bu\in\HH$ a ground state
solution of the problem
\begin{equation*}
\left\{
\begin{array}{ll}
-\Delta u_1 +V_{1,\infty} u_1 - |u_1|^{p-1}u_1=\g u_2 u_3	& \hbox{in }\RN,
\\
-\Delta u_2 +V_{2,\infty} u_2 - |u_2|^{p-1}u_2=\g u_1 u_3	& \hbox{in }\RN,
\\
-\Delta u_3 +V_{3,\infty} u_3 - |u_3|^{p-1}u_3=\g u_1 u_2	& \hbox{in }\RN,
\\
u_1,u_2,u_3\in \H.
\end{array}
\right.
\end{equation*}
Let $t_{\bu}>0$ be such that $t_{\bu}\bu\in \Ne_{\bf V}$. By ({\bf V2}), we have
\begin{align*}
c_\infty & =I_\infty (\bu) \ge I_\infty \big(t_{\bu}\bu\big)
\\
&=I_{\bf V}\big(t_{\bu}\bu\big) + t_{\bu}^2 \SU \irn \big(V_{i,\infty} -V_i(x)\big) u_i^2
> c_{\bf V},
\end{align*}
and then we conclude.
\end{proof}

\subsection{Proof of Theorem \ref{main3}}

Let $\{\bun\}_n \subset \Ne_{\bf V}$ such that
\begin{equation}\label{eq:lim2}
\lim_n I_{\bf V} (\bun)=c_{\bf V}.
\end{equation}
We define the functional $J\colon \HH\to\R$ as:
\begin{equation*}
J(\bu)= \SU  \irn \frac {1}{6}\big(|\n u_i|^2 +  V_i(x) u_i^2\big) + \frac{p-2}{3(p+1)} |u_i|^{p+1}.
\end{equation*}
Observe that for any $\bu\in\Ne_{\bf V},$ we have $I_{\bf V} (\bu) =J(\bu).$
\\
To prove Theorem \ref{main3}, we need some compactness properties on the
sequence $\{\bun\}_n.$\\
We denote by $\nu_n$ the measure
\begin{equation}\label{eq:meas}
\nu_n(\O)=  \SU  \int_\O \frac {1}{6}\big(|\n u_i|^2 +  V_i(x) u_i^2\big) +  \frac{p-2}{3(p+1)}|u_i|^{p+1}.
\end{equation}
By \eqref{eq:lim2} we have
\[
\nu_n(\RN)=J(\bun)\to c_{\bf V}
\]
and then, by P.L.~Lions \cite{L1}, there are three possibilities:
\begin{description}
\item[\;\;{\it vanishing}\,{\rm :}] for all $r>0$
\[
\lim_n \sup_{\xi \in \RN}\int_{B_r(\xi)} d \nu_n =0;
\]
\item[\;\;{\it dichotomy}\,{\rm :}] there exist a constant $\tilde
c\in (0, c_{\bf V})$, two sequences $\{\xi_n\}_n$ and $\{r_n\}_n$, with $r_n
\to +\infty$ and two nonnegative measures $\nu_n^1$ and $\nu_n^2$
such that
\begin{align*}
0\le \nu_n^1 + \nu_n^2 \le \nu_n,&\qquad \nu_n^1(\RN) \to \tilde
c,\;\; \nu_n^2(\RN) \to c_{\bf V} -\tilde c,
\\
\hbox{supp}(\nu_n^1)\subset B_{r_n}(\xi_n),&\qquad
\hbox{supp}(\nu_n^2)\subset \RN \setminus B_{2r_n}(\xi_n);
\end{align*}
\item[\;\;{\it compactness}\,{\rm :}] there exists a sequence
$\{\xi_n\}_n$ in $\RN$ with the following property: for any $\d>0$,
there exists $r=r(\d)>0$ such that
\[
\int_{B_r(\xi_n)} d \nu_n \ge c_{\bf V} -\d.
\]
\end{description}
Arguing as in \cite{WZ}, we prove the following
\begin{lemma}\label{le:concentr}
Compactness holds
for the sequence of measures $\{\nu_n\}_n$, defined in \eqref{eq:meas}.
\end{lemma}
\begin{proof}
{\sc Vanishing does not occur}
\\
Suppose by contradiction, that for all $r>0$
\[
\lim_n \sup_{\xi \in \RN}\int_{B_r(\xi)} d \nu_n =0.
\]
In particular, we deduce that there exists $\bar r>0$ such that
\begin{equation*}
\lim_n \sup_{\xi \in \RN}\SU \int_{B_{\bar r}(\xi)} u_{i,n}^2=0.
\end{equation*}
Arguing as in Lemma \ref{le:nonvan}, we get a contradiction.
\\
\\
{\sc Dichotomy does not occur}
\\
Suppose by contradiction that there exist a constant $\tilde c\in
(0, c_{\bf V})$, two sequences $\{\xi_n\}_n$ and $\{r_n\}_n$, with $r_n \to
+\infty$ and two nonnegative measures $\nu_n^1$ and $\nu_n^2$ such
that
\begin{align*}
0\le \nu_n^1 + \nu_n^2 \le \nu_n,&\qquad \nu_n^1(\RN) \to \tilde
c,\;\; \nu_n^2(\RN) \to c_{\bf V} -\tilde c,
\\
\hbox{supp}(\nu_n^1)\subset B_{r_n}(\xi_n),&\qquad
\hbox{supp}(\nu_n^2)\subset \RN \setminus B_{2r_n}(\xi_n).
\end{align*}
Let $\rho_n \in C^1(\RN)$ be such that $\rho_n\equiv 1$ in $B_{r_n}(\xi_n)$, $\rho_n\equiv 0$ in $\RN \setminus
B_{2r_n}(\xi_n)$, $0\le \rho_n\le 1$ and $|\n \rho_n|\le 2/r_n$.
Set
$v_{i,n}=\rho u_{i,n}$ and $w_{i,n}=(1-\rho)u_{i,n}$, and we denote $\bvn=(v_{1,n},v_{2,n},v_{3,n})$ and $\bwn=(w_{1,n},w_{2,n},w_{3,n})$.
\\
It is easy to see that
\begin{align*}
\liminf_n J(\bvn) &\ge \tilde c,
\\
\liminf_n J(\bwn) &\ge c_{\bf V} - \tilde c.
\end{align*}
Moreover, denoting $\O_n:= B_{2r_n}(\xi_n)\setminus
B_{r_n}(\xi_n)$, we have
\begin{equation*}
\nu_n(\O_n) \to 0, \qquad \hbox{as }n\to \infty,
\end{equation*}
namely
\begin{align}
\int_{\O_n}|\n u_{i,n}|^2 +   V_i(x)u_{i,n}^2 \to 0,&\quad\hbox{ for all }i=1,2,3, \label{eq:1}
\\
\int_{\O_n} |u_{i,n}|^{p+1} \to 0 ,&\quad\hbox{ for all } i=1,2,3.
\label{eq:2}
\end{align}
Moreover by \eqref{eq:2}, we have
\begin{equation}\label{eq:3}
\int_{\O_n} \uun\to 0.
\end{equation}
By simple computations, by \eqref{eq:1}, \eqref{eq:2} and \eqref{eq:3} we
infer that for all $i=1,2,3$
\begin{align*}
&\int_{\O_n}|\n v_{i,n}|^2 +   V_i(x) v^2_{i,n} \to 0 , 
&\int_{\O_n}&|\n w_{i,n}|^2 +   V_i(x) w^2_{i,n} \to 0 ,
\\
&\int_{\O_n}| v_{i,n}|^{p+1}\to 0 , &\int_{\O_n}&| w_{i,n}|^{p+1} \to 0 ,
\\
&\int_{\O_n} v_{1,n}v_{2,n}v_{3,n}\to 0 , &\int_{\O_n}& w_{1,n}w_{2,n}w_{3,n}\to 0 .
\end{align*}
Hence, we deduce that  for all $i=1,2,3$
\begin{align}
\irn |\n u_{i,n}|^2 +   V_i(x) u^2_{i,n} &=\irn|\n v_{i,n}|^2 + V_i(x)v^2_{i,n}
\nonumber
\\
&\quad+\irn|\n w_{i,n}|^2 + V_i(x) w^2_{i,n} +o_n(1), \label{eq:normV}
\\
\irn |u_{i,n}|^{p+1}&=\irn|v_{i,n}|^{p+1}+\irn|w_{i,n}|^{p+1}
+o_n(1)\label{eq:p+1V},
\\
\irn \uun&=\irn v_{1,n}v_{2,n}v_{3,n}+ \irn w_{1,n}w_{2,n}w_{3,n}+o_n(1). \label{eq:uuuV}
\end{align}
Hence, by \eqref{eq:normV} and \eqref{eq:p+1V}, we get
\begin{align}
J(\bun)= J(\bvn)+J(\bwn)+o_n(1).\nonumber
\end{align}
Then
\begin{align*}
c_{\bf V}=\lim_n J(\bun)
\ge\liminf_n J(\bvn)
+\liminf_n J(\bwn)
\ge \tilde c+(c_{\bf V}-\tilde c)=c_{\bf V},
\end{align*}
hence
\begin{align}
\lim_n J(\bvn)&=\tilde c, \label{eq:jv}
\\
\lim_n J(\bwn)&=c_{\bf V}-\tilde c. \nonumber
\end{align}
Let us observe, moreover, that by
\eqref{eq:normV}, \eqref{eq:p+1V}  and \eqref{eq:uuuV}, we have
\begin{equation}
0=G_{\bf V}(\bun) = G_{\bf V}(\bvn)+G_{\bf V}(\bwn)+o_n(1)\label{eq:GV}.
\end{equation}
We have to distinguish three cases.
\\
\
\\
{\sc Case 1:} up to a subsequence, $G_{\bf V}(\bvn) \le 0$.
\\
By Lemma \ref{le:NV3}, for any $n \ge 1$, there exists $\t_n>0$ such
that $\t_n \bvn \in \Ne_{\bf V}$, and then
\begin{equation}\label{eq:nehariV}
\SU\irn |\n v_{i,n}|^2 + V_i(x) v_{i,n}^2 - \t_n^{p-1} |v_{i,n}|^{p+1}= 3\g\t_n\irn v_{1,n}v_{2,n}v_{3,n}.
\end{equation}
By \eqref{eq:nehariV} we have
\begin{align*}
\SU \irn (\t_n-1)\big(|\n v_{i,n}|^2 + V_i(x) v_{i,n}^2\big) +
(\t_n^{p-1}-\t_n)|v_{i,k}|^{p+1} \le 0,
\end{align*}
and, by ({\bf V3}), we deduce that $\t_n\le 1$. Therefore, for all
$n\ge 1,$ by ({\bf V3}) and  \eqref{eq:jv},
\begin{align*}
c_{\bf V} & \le  I_{\bf V}( \t_n \bvn) = J( \t_n \bvn) \le J(\bvn) \to \tilde c<c_{\bf V},
\end{align*}
which is a contradiction.
\\
\
\\
{\sc Case 2:} up to a subsequence, $ G_{\bf V}(\bwn) \le 0$.
\\
We can argue as in the previous case.
\\
\
\\
{\sc Case 3:} up to a subsequence, $G_{\bf V}(\bvn) > 0$ and $G_{\bf V}(\bwn)> 0$.
\\
By \eqref{eq:GV}, we infer that $G_{\bf V}(\bvn)=o_n(1)$ and $G_{\bf V}(\bwn)=o_n(1)$. 
\\
Since $\{\bun\}_n\subset \Ne_{\bf V}$, by \eqref{eq:lim2} and by Lemma \ref{le:>d}, we infer that there exists $\mu>0$ and $n_0\in \N$ such that
\[
\SU \irn |u_{i,n}|^{p+1}>\mu, \hbox{ for }n\ge n_0,
\]
and hence, by \eqref{eq:p+1V}, there exists $n_1\in \N$ such that
\[
\SU \irn |v_{i,n}|^{p+1}>\frac \mu 2, \hbox{ or }\SU \irn |w_{i,n}|^{p+1}>\frac \mu 2, \hbox{ for }n\ge n_1.
\]
Without loss of generality, suppose that 
\begin{equation} \label{eq:mu}
\SU \irn |v_{i,n}|^{p+1}>\frac \mu 2, \hbox{ for }n\ge 1.
\end{equation}
Let $\{\t_n\}_n$ be such that $\t_n \bvn \in \Ne_{\bf V}$. Combining \eqref{eq:normV} and \eqref{eq:mu}, we deduce that $\{\t_n\}_n$ is bounded. If $\t_n\le 1+o_n(1)$, we can repeat the arguments of Case 1. Suppose, therefore, that
\begin{equation*}
\lim_n \t_n=\t_0>1.
\end{equation*}
We have
\begin{align*}
o_n(1)&= G_{\bf V}(\bvn) \!= \!\SU\irn \!|\n v_{i,n}|^2+  V_i(x) v_{i,n}^2- |v_{i,n}|^{p+1}
-3\g \! \irn \!\! v_{1,n} v_{2,n} v_{3,n}
\\
&=\left( 1-\frac{1}{\t_n}\right)\SU\irn |\n v_{i,n}|^2+ V_i(x) v_{i,n}^2 
+\left(\t_n^{p-2}- 1\right) \SU\irn |v_{i,n}|^{p+1}
\end{align*}
and so
\[
\irn |\n v_{i,n}|^2 +  V_i(x) v_{i,n}^2=o_n(1),\hbox{ for }i=1,2,3,
\]
which contradicts \eqref{eq:mu}.
\end{proof}

\begin{proofmain3}
Let $\{\bun\}_n$ be a sequence in $\Ne_{\bf V}$ such that \eqref{eq:lim2} holds. 
\\
By ({\bf V3}) and \eqref{eq:lim2}, we deduce that $\{\bun\}_n$ is
bounded in $\HH,$ so there exists $(\bar u_1,\bar u_2,\bar u_3)\in\HH$ such that, up to a
subsequence,
\begin{align}
&u_{i,n}\rightharpoonup \bar u_i\quad\hbox{weakly in }\H,
\label{eq:weakmain2}
\\
&u_{i,n}\to \bar u_i\quad\hbox{in }L^s_{loc}(\RN), \hbox{ with }1\le s<2^*.  \label{eq:loc}
\end{align}
We define the measures $\{\nu_n\}_n$ as in \eqref{eq:meas}; by Lemma \ref{le:concentr}
there exists a sequence $\{\xi_n\}_n$
in $\RN$ with the following property: for any $\d>0$, there exists
$r=r(\d)>0$ such that
\begin{equation}\label{eq:bc}
\SU  \int_{B_r^c(\xi_n)} \frac {1}{6}\big(|\n u_i|^2 +  V_i(x) u_i^2\big) +  \frac{p-2}{3(p+1)}|u_i|^{p+1}< \d.
\end{equation}
\noindent{\sc Claim}: $\{\xi_n\}_n$ is bounded in $\RN$.
\\
Suppose by contradiction that, up to a subsequence, $|\xi_n|\to \infty$, as $n\to \infty$.
\\
Fix $\mu>0$ and let ${\bf V}_{\infty,\mu}=(V_{1,\infty}-\mu,V_{2,\infty}-\mu,V_{3,\infty}-\mu)$ and $I_{\infty,\mu}$ the functional associated with ${\bf V}_{\infty,\mu}$. For any $n\ge1,$ let $\bzn=(z_{1,n},z_{2,n},z_{3,n})$, where $z_{i,n}=u_{i,n}(\cdot - \xi_n)$, and $t_n>0$ such that the functions $t_n \bzn$ are in the Nehari manifold of $I_{\infty,\mu}$.
\\
Let $\d>0$ and consider $r>0$ such that \eqref{eq:bc} holds. For $n$ sufficiently large, we have
\[
V_i(x+\xi_n)\ge V_{i,\infty}-\mu, \qquad \hbox{for all }x\in B_r.
\]
Hence we have
\begin{align*}
c_{\bf V}+o_n(1) &= I_{\bf V}(\bun) \ge I_{\bf V}(t_n \bun) 
\\
&= I_{\infty,\mu}(t_n \bun) +
\frac{t_n^2}{2}\SU \irn \left(V_i(x)- (V_{i,\infty}-\mu)\right)u_{i,n}^2
\\
&\ge c_{{\bf V}_{\infty,\mu}} +
\frac{t_n^2}{2}\SU \int_{B_r} \left(V_i(x+\xi_n)- (V_{i,\infty}-\mu)\right)z_{i,n}^2
\\
&\quad+ \frac{t_n^2}{2}\SU\int_{B_r^c} \left(V_i(x+\xi_n)- (V_{i,\infty}-\mu)\right)z_{i,n}^2
\\
&\ge c_{{\bf V}_{\infty,\mu}} 
- \frac{t_n^2}{2}\SU\int_{B_r^c} \left|V_i(x+\xi_n)- (V_{i,\infty}-\mu)\right|z_{i,n}^2.
\end{align*}
Since by \eqref{eq:bc}
\[
\SU\int_{B_r^c} \left|V_i(x+\xi_n)- (V_{i,\infty}-\mu)\right|z_{i,n}^2\le
C\d,\quad\hbox{for any }n\ge 1,
\]
and $\{t_n\}_n$ is bounded (the proof is the same as in Lemma \ref{le:t-bdd}), we get that $c_{\bf V} \ge c_{{\bf V}_{\infty,\mu}}  - C\d$. By the
arbitrariness in the choice of $\d>0,$ we have $c_{\bf V} \ge c_{{\bf V}_{\infty,\mu}} .$
Using Lemma \ref{le:vn} we conclude that $c_{\bf V}\ge c_\infty,$ which
contradicts Lemma~\ref{le:cinfty}.
\\
\
\\
So $\{\xi_n\}_n$ is bounded in $\RN$ and then, by \eqref{eq:bc}, for any
$\delta>0$ there exists $r>0$ such that
\begin{equation}\label{eq:conc2}
\|u_{i,n}\|_{H^1(B_r^c)}<\d,\quad \hbox{uniformly for }n\ge 1, \quad i=1,2,3.
\end{equation}
By \eqref{eq:weakmain2}, \eqref{eq:loc} and \eqref{eq:conc2}, we have that, taken $s\in [2,2^*[$,
for any $\d>0$ there exists $r>0$ such that, for any $n\ge 1$
large enough
\begin{align*}
\|u_{i,n}-\bar u_i\|_{L^s(\RN)}
& \le \|u_{i,n}-\bar u_i\|_{L^s(B_r)}+\|u_n-\bar u_i\|_{L^s(B^c_r)}
\\
& \le \d + C\left( \|u_{i,n}\|_{H^{1}(B_r^c)}+\|\bar u_i\|_{H^1(B_r^c)}\right)\le
(1+2C)\d,
\end{align*}
where $C>0$ is the constant of the embedding
$H^{1}(B_r^c)\hookrightarrow L^s(B^c_r).$
We deduce that
\begin{equation}\label{eq:conv}
u_{i,n}\to \bar u_i\hbox{ in }L^s(\RN),\;\hbox{for any }s\in [2,2^*[.
\end{equation}
By \eqref{eq:lim2}, we can suppose (see \cite{Wi}) that $\{\bun\}_n$
is a Palais-Smale sequence for ${I_{\bf V}}_{|\Ne_{\bf V}}$ and, as a consequence,
it is easy to see that $\{\bun\}_n$ is a Palais-Smale sequence for
$I_{\bf V}$. By \eqref{eq:weakmain2} and \eqref{eq:conv},
we conclude that $I'_{\bf V}(\bar \bu)=0.$
\\
We have to show that $\bar \bu\neq (0,0,0)$. Suppose by contradiction that $\bar \bu=(0,0,0)$, since $\{\bun\}_n$ is in $\Ne_{\bf V}$ and by \eqref{eq:conv}, we should have that $\|\bun\|\to 0$ and this is impossible by Lemma \ref{le:NV1}.
Hence $\bar \bu\neq (0,0,0)$ and so $\bar \bu\in \Ne_{\bf V}$.
\\
Finally, by \eqref{eq:lim2}, \eqref{eq:weakmain2} and \eqref{eq:conv} and by ({\bf V2-3}) we get
\begin{align*}
c_{\bf V} \le I_{\bf V}(\bar \bu) \le \liminf I_{\bf V}(\bun) = c_{\bf V},
\end{align*}
so we can conclude that $\bar \bu$ is a ground state solution of \eqref{eqv}.
\end{proofmain3}

\subsection{Proof of Theorem \ref{main4}}

In this section, we assume that ${\bf V}$ satisfies ({\bf V2'}) instead of ({\bf V2}). We start proving the analogous of Lemma \ref{le:cinfty}.
\begin{lemma}\label{le:cinfty2}
If ${\bf V}$ satisfies ({\bf V1}), ({\bf V2'}) and ({\bf V3}), then there exists $\g_0>0$ such that, for any $\g\in \R$ with $|\g|>\g_0$, we get $c_{\bf V}<c_\infty$.
\end{lemma}

\begin{proof}
By Theorem \ref{main2}, exists $\g_0>0$ such that, for any $\g\in \R$ with $|\g|>\g_0$, there exists $\bu=\uuu\in\HH$ a vector ground state
solution of the problem
\begin{equation*}
\left\{
\begin{array}{ll}
-\Delta u_1 +V_{1,\infty} u_1 - |u_1|^{p-1}u_1=\g u_2 u_3	& \hbox{in }\RN,
\\
-\Delta u_2 +V_{2,\infty} u_2 - |u_2|^{p-1}u_2=\g u_1 u_3	& \hbox{in }\RN,
\\
-\Delta u_3 +V_{3,\infty} u_3 - |u_3|^{p-1}u_3=\g u_1 u_2	& \hbox{in }\RN,
\\
u_1,u_2,u_3\in \H.
\end{array}
\right.
\end{equation*}
Let $t_{\bu}>0$ be such that $t_{\bu}\bu\in \Ne_{\bf V}$. By ({\bf V2'}) and since $\bu$ is a vector ground state and so $u_i\neq 0$, for all $i=1,2,3$, we have
\begin{align*}
c_\infty & =I_\infty (\bu) \ge I_\infty \big(t_{\bu}\bu\big)
\\
&=I_{\bf V}\big(t_{\bu}\bu\big) + t_{\bu}^2 \SU \irn \big(V_{i,\infty} -V_i(x)\big) u_i^2
> c_{\bf V},
\end{align*}
and then we conclude.
\end{proof}

Let us now prove the existence of vector ground state.

\begin{proofmain4}
By Lemma \ref{le:cinfty2}, we infer that for $|\g|$ sufficiently large $c_{\bf V}<c_\infty$. We can repeat the arguments of the proof of Theorem \ref{main3} to prove the existence of a ground state $\bu$ of \eqref{eqv}. Therefore, we have only to show that $\bu$ is a vector solution. 
\\
By \cite{R}, there exists a ground state, $\bar u_i\in \H$, of 
\[
-\Delta u+V_i(x) u= |u|^{p-1}u, \quad \hbox{in }\RN, \qquad u \in \H,
\]
and moreover, it is easy to see that $\bar u_i>0$. 
\\
Now the proof is similar to that of Theorem \ref{main2} and so we omit it.
\end{proofmain4}

\end{document}